\documentclass[12pt]{article}

\usepackage{amsmath}
\usepackage{amssymb}
\usepackage{microtype}
\usepackage[margin=1.30in]{geometry}
\usepackage{amsthm}
\usepackage{url}

\newtheorem{thm}{Theorem}[section]
\newtheorem{defn}[thm]{Definition}
\newtheorem{lemma}[thm]{Lemma}

\newtheorem{cor}[thm]{Corollary}

\newtheorem{conj}[thm]{Conjecture}

\newcommand{\cE}[0]{{\cal E}}

\newcommand{\cI}[0]{{\cal I}}

\newcommand{\cO}[0]{{\cal O}}

\newcommand{\ga}[0]{\alpha}

\newcommand{\gl}[0]{\lambda}

\newcommand{\gld}[0]{t\left(1-\frac{t}{2^{d-1}}\right)^{d-1}}

\begin{document}
\renewcommand{\thefootnote}{\fnsymbol{footnote}}

\title{The independent set sequence of regular bipartite graphs}

\author{David Galvin\thanks{Department of Mathematics,
University of Notre Dame, 255 Hurley Hall, Notre Dame IN
46556; dgalvin1@nd.edu. Research supported in part by National Security Agency grant H98230-10-1-0364.}}

\maketitle

\begin{abstract}

Let $i_t(G)$ be the number of independent sets of size $t$ in a graph $G$. Alavi, Erd\H{o}s, Malde and Schwenk made the conjecture that if $G$ is a tree then the independent set sequence $\{i_t(G)\}_{t\geq 0}$ of $G$ is unimodal; Levit and Mandrescu further conjectured that this should hold for all bipartite $G$.

We consider the independent set sequence of finite regular bipartite graphs, and graphs obtained from these by percolation (independent deletion of edges).
Using bounds on the independent set polynomial $P(G,\lambda):=\sum_{t \geq 0} i_t(G)\lambda^t$ for these graphs, we obtain partial unimodality results in these cases.

We then focus on the discrete hypercube $Q_d$, the graph on vertex set $\{0,1\}^d$ with two strings adjacent if they differ on exactly one coordinate. We obtain asymptotically tight estimates for $i_{t(d)}(Q_d)$ in the range
$t(d)/2^{d-1} > 1-1/\sqrt{2}$, and nearly matching upper and lower bounds otherwise. We use these estimates to obtain a stronger partial unimodality result for the independent set sequence of $Q_d$.

\end{abstract}


\section{Introduction and statement of results}

For a (finite, simple, undirected, loopless) graph $G=(V,E)$ set
$$
i_t(G) =\{I \in \cI(G) : |I|=t\}
$$
where $\cI(G)$ is the collection of independent sets of $G$ (sets of vertices spanning no edges). The {\em independent set sequence} of $G$ is the sequence $\{i_t(G)\}_{t=0}^{\alpha(G)}$ where $\alpha(G)$ is the size of the largest independent set in $G$. The {\em independent set polynomial} or {\em stable set polynomial} of $G$, first introduced explicitly by Gutman and Harary \cite{GutmanHarary}, is the polynomial
$$
P(G, \lambda) = \sum_{t=0}^{\alpha(G)} i_t(G)\lambda^t.
$$

A sequence $\{a_i\}_{i=0}^n$ is said to be {\em unimodal} (with {\em mode at} $k$) if
$$
a_0 \leq a_1 \leq \ldots \leq a_k \geq a_{k+1} \geq \ldots \geq a_n,
$$
and the polynomial $\sum_{i=0}^n a_i\lambda^i$ is said to be unimodal if its sequence of coefficients is.
There has been some consideration in the literature of the question of the unimodality of the independent set polynomial of a graph. There are two major positive results.
The first of these follows from the celebrated result of Heilmann and Lieb \cite{HeilmannLieb} to the effect that for any graph $G$ the matching polynomial $\sum_{i\geq 0} m_t(G)\lambda^t$ of $G$ (where $m_t(G)$ is the number of matchings of $G$ of size $t$) has only real roots, implying by a theorem of Newton (see for example \cite[page 504]{Stanley}) that it is unimodal
(in fact, log-concave). Since independent sets of a fixed size in the line graph of a graph are in bijection with matchings of that size in the original graph, Heilmann and Lieb's result shows that if $G$ is the line graph of a graph, then $P(G,\lambda)$ is unimodal.

The second of these positive results is due to
Hamidoune \cite{Hamidoune}, who showed that if $G$ is claw-free (that is, does not contain a star on four vertices as an induced subgraph) then $P(G,\lambda)$ is unimodal. Later Chudnovsky and Seymour \cite{ChudnovskySeymour} showed that in this case  $P(G,\lambda)$ also has all real roots. Since line graphs are claw-free, this result generalizes that of Heilmann and Lieb.

On the other hand, if $G$ contains a claw then $P(G,\lambda)$ may not be unimodal. For example, the graph obtained from a claw by replacing each of the degree $1$ vertices with a $K_4$ (the complete graph on $4$ vertices), the degree $3$ vertex with a $K_{37}$, and all edges with complete bipartite graphs, has independent set polynomial $1+49\lambda + 48\lambda^2 +
64\lambda^3$. In fact, Alavi, Erd\H{o}s, Malde and Schwenk \cite{AlaviErdosMaldeSchwenk} showed that in general the independent set polynomial of a graph can have every possible pattern of increases and decreases. Specifically, they showed that for any integer $m \geq 1$ and any permutation $\pi$ of $\{1, \ldots, m\}$, there is a graph $G$ with $\alpha(G)=m$ and with
\begin{equation} \label{seq-AEMS}
i_{\pi(1)}(G) < i_{\pi(2)}(G) < \ldots < i_{\pi(m)}(G).
\end{equation}
In the language of \cite{AlaviErdosMaldeSchwenk}, the independent set polynomial of a general graph is {\em unconstrained}. Note that there are at most $2^m=o(m!)$ permutations $\pi$ of $\{1,\ldots, m\}$ for which (\ref{seq-AEMS}) holds for some graph $G$ satisfying $\alpha(G)=m$ and $P(G,\lambda)$ unimodal.

Alavi, Erd\H{o}s, Malde and Schwenk made the following positive conjecture.
\begin{conj} \label{conj-fromAEMS}
If $T$ is a tree then $P(T,\lambda)$ is unimodal.
\end{conj}
In \cite{LevitMandrescu} Levit and Mandrescu make the stronger conjecture that $P(G, \gl)$ is unimodal if $G$ is a K\"{o}nig-Egerv\'{a}ry graph (a graph in which the size of the largest independent set plus the size of the largest matching equals the number of vertices in the graph). In particular, since all bipartite graphs are K\"{o}nig-Egerv\'{a}ry graphs, we have the following:
\begin{conj} \label{conj-bip}
If $G$ is bipartite then $P(G,\lambda)$ is unimodal.
\end{conj}
Very little progress has been made towards this conjecture; indeed, even Conjecture \ref{conj-fromAEMS} remains open. A partial result of Levit and Mandrescu \cite{LevitMandrescu} is that for $G$ a bipartite graph, the final third of the coefficients of $P(G,\gl)$ form a decreasing sequence, that is,
\begin{equation} \label{seq-LM}
i_{\lceil (2\ga(G)-1)/3\rceil}(G) \geq i_{\lceil (2\ga(G)-1)/3\rceil+1}(G) \geq \ldots \geq i_{\ga(G)}(G).
\end{equation}

\medskip

In this note we consider the independent set polynomial of a graph drawn from the families
of regular and almost-regular bipartite graphs. One approach to showing unimodality of the independent set polynomial of a graph is to obtain upper and lower bounds on $i_t(G)$ for each $0 \leq t \leq \alpha(G)$, and this is the approach that we take here. We begin with two simple bounds. Here and throughout $H(x)=-x\log x -(1-x)\log (1-x)$ is the binary entropy function and $\log=\log_2$.

\begin{lemma} \label{lem-fixed_count_u_l_bounds}
For $d$-regular $G$ (not necessarily bipartite),
\begin{equation} \label{eq-fromCGT}
i_t(G) \leq \exp_2\left\{H\left(\frac{2t}{|V|}\right)\frac{|V|}{2}+ \frac{|V|}{2d}\right\}.
\end{equation}
If $G$ is bipartite then we also have
\begin{equation} \label{eq-sillylb}
i_t(G) \geq {\frac{|V|}{2} \choose t} \geq \exp_2\left\{H\left(\frac{2t}{|V|}\right)\frac{|V|}{2} -\frac{1}{2}\log|V|\right\}
\end{equation}
\end{lemma}

\medskip

\noindent {\em Proof}:
We begin with (\ref{eq-fromCGT}), which is based on the observation that for all $t>0$ and $\lambda>0$ we have $i_t(G) \lambda^t \leq P(G, \lambda)$ and so
\begin{equation} \label{ins1}
i_t(G) \leq \min_{\lambda>0} \left\{\frac{P(G,\lambda)}{\lambda^t}\right\}.
\end{equation}
We now use the following inequality, proved in \cite{Zhao}: for $d$-regular $G$ (not necessarily bipartite)
\begin{equation} \label{eq-general_ub_on_part}
P(G, \lambda) \leq 2^\frac{|V|}{2d}(1+\lambda)^{\frac{|V|}{2}}.
\end{equation}
(A weaker bound with $|V|/d$ replacing $|V|/(2d)$ had earlier been obtained in \cite{CarrollGalvinTetali}). Taking $\lambda = \frac{2t}{|V(G)|-2t}$ in (\ref{eq-general_ub_on_part}) and plugging into (\ref{ins1}) we get (\ref{eq-fromCGT}) for $t \neq 0, |V|/2$, and these two extreme cases are trivial.

The near-matching lower bound (\ref{eq-sillylb}) is obtained by specifying a bipartition $V = {\mathcal E} \cup {\mathcal O}$ of $G$ and considering only those independent sets which are subsets of ${\mathcal E}$. The second inequality in (\ref{eq-sillylb}) follows for all $|V|\geq 1$ and $0 \leq t \leq |V|/2$ from Stirling's formula. (Note that for regular bipartite $G$, $\alpha(G)=|V(G)|/2$.)
\qed

\medskip

The upper bound (\ref{eq-fromCGT}) is close to best possible; the graph consisting of disjoint copies of $K_{d,d}$, the complete bipartite graph with $d$ vertices in each class, shows that $|V|/2d$ cannot be replaced by $c|V|/d$ for any $c < 1/2$.

Combining (\ref{eq-fromCGT}) and (\ref{eq-sillylb}) we see that for regular, bipartite $G$ we have $i_t(G) \approx {|V|/2 \choose t}$. This suggests that in contrast to the general situation, for regular bipartite $G$ the independent set polynomial may be quite constrained (in the language of \cite{AlaviErdosMaldeSchwenk}), in the sense that any sufficiently sparse subsequence of $\{i_t(G)\}_{t=0}^{|V(G)|/2}$ should be unimodal.


To state a precise result in this direction, it will be helpful to set up some notation.
\begin{defn}
A sequence $\{a_i\}_{i\geq 0}$ is {\em $s$-step monotone increasing} on the interval $[a,b]$ (where $a$ and $b$ satisfy $0 \leq a \leq b$) if for every $a \leq i \leq j \leq b$ with $j-i \geq s$ we have $a_i \leq a_j$. We define {\em $s$-step monotone decreasing} analogously.
\end{defn}
\noindent Note that $\{a_i\}_{i=0}^r$ being unimodal with mode at $k$ is equivalent to $\{a_i\}_{i=0}^r$ being $1$-step monotone increasing on $[0,k]$ and $1$-step monotone decreasing on $[k,r]$.
\begin{defn}
A bipartite graph $G$ on $2n$ vertices has {\em property $(\beta,\gamma,s)$} if the sequence $\{i_t(G)\}_{t \geq 0}$ is $s$-step monotone increasing on the interval $[\beta n, (1-\gamma)n/2]$ and $s$-step monotone decreasing on the interval $[(1+\gamma)n/2, (1-\beta)n]$.
\end{defn}
\noindent Note that property $(0,0,1)$ is equivalent to the unimodality of $\{i_t(G)\}_{t = 0}^n$ with mode at $n/2$.


\medskip

The first aim of this paper is to prove a number of results establishing property $(\beta,\gamma,s)$ for $\beta$ and $\gamma$ arbitrarily small constants and $s=o(n)$.
\begin{thm} \label{thm-step_unimodality}
Fix $\varepsilon > 0$. There is $C=C(\varepsilon)>0$ such that if $G$ is a $2n$-vertex, $d$-regular bipartite graph, then $G$ has property $(0,\varepsilon,s)$ with
$$
s = C\max\left\{\log n,\frac{n}{d}\right\}.
$$
\end{thm}
\noindent We prove Theorem \ref{thm-step_unimodality} in Section \ref{sec-perc}. Note that if $d=\omega(1)$ above, then $s=o(n)$.

The condition that $G$ be regular can be relaxed quite a bit; using a recent result in the spirit of (\ref{eq-general_ub_on_part}) obtained in \cite{EngbersGalvin} we may extend Theorem \ref{thm-step_unimodality} to all bipartite graphs which are suitably close to regular. Specifically, let $G$ be a bipartite graph on $2n$ vertices with bipartition classes ${\mathcal E}$ and ${\mathcal O}$ (with $|\cO|\geq |\cE|$), let $d$ be an arbitrary positive parameter and let
$$
h(G,d) = \frac{1}{d} + \frac{|\{v \in \cE: d(v) < d\}|}{n}  +  \frac{1}{dn}\sum_{v \in \cO} (d(v)-d)\mathbf{1}_{\{d(v)\geq d\}} + \frac{|\cO|-|\cE|}{n}
$$
(where $\mathbf{1}_E$ is the indicator function of the event $E$). The parameter $h(G,d)$ is intended to capture the extent
to which $G$ is ``almost $d$-regular'': $h(G,d)$ being close to $1/d$ (the value it takes when $G$ is $d$-regular) means that $G$ has not too many low degree vertices, that the sum of the
degrees of high degree vertices is not too large, and that the difference between
the sizes of the partition classes is not too great.

Similar to Theorem \ref{thm-step_unimodality}, we have the following statement, whose proof is also given in Section \ref{sec-perc}.
\begin{thm} \label{thm-step_unimodality-gen}
Fix $\varepsilon > 0$. There is $C=C(\varepsilon)>0$ such that if $G$ is a $2n$-vertex bipartite graph and $d$ is an arbitrary positive parameter, then $G$ has property $(\varepsilon, \varepsilon, s)$ with
$$
s = C\max\left\{\log n,n h(G,d)\right\}.
$$
\end{thm}
This result is mainly of interest when $s=o(n)$, that is, when there is a choice of $d$ for which $h(G,d)=o(1)$. One quite general situation in which this occurs is in percolation. Given a graph $G$ and a parameter $0 \leq p \leq 1$, let $G^p$ be a random subgraph of $G$ obtained by deleting each edge independently with probability $1-p$ (so
$$
\Pr\left(G^p=H\right) = p^{|E(H)|}(1-p)^{|E(G)|-|E(H)|}).
$$
In \cite[Section 4]{EngbersGalvin} it is shown that if $G$ is a $d$-regular bipartite graph on $2n$ vertices, and $G^p$ is obtained from $G$ by percolation with $p \geq f(d)/d$ for an arbitrary function $f(d)=\omega(1)$, then there is a function $g(d)=o(1)$ such that with probability at least $1-g(d)$, we have that $h(G^p,d')\leq g(d)$ (where $d'=dp-(2dp)^{1/2}f(d)^{1/4}$). We thus have the following corollary of Theorem \ref{thm-step_unimodality-gen}.
\begin{cor} \label{cor-perc}
Fix $\varepsilon > 0$. Let $G$ be a $d$-regular bipartite graph on $2n$ vertices. Let $p=f(d)/d$ for an arbitrary function $f(d)=\omega(1)$. There are functions $g(d)=o(1)$ and $s(n,d)=o(n)$ such that with probability at least $1-g(d)$, $G^p$ has property $(\varepsilon, \varepsilon, s(n,d))$.
\end{cor}
A particularly interesting application of Corollary \ref{cor-perc} is to the random equi-bipartite graph $G(n,n,p)$. This is the graph on vertex set ${\mathcal E} \cup {\mathcal O}$ with $|{\mathcal E}|=|{\mathcal O}|=n$ in which the edge $\{u,v\}$ ($u \in {\mathcal E}$, $v \in {\mathcal O}$) is present with probability $p$, independently for all choices of $u$ and $v$. (In other words, $G(n,n,p)$ is the result of percolation on $K_{n,n}$ with parameter $p$.)
\begin{cor} \label{cor-random-bip-graph}
Fix $\varepsilon > 0$. Let $p=\omega(1)/n$. There is a function $s(n)=o(n)$ so that almost surely (with probability tending to $1$ as $n$ goes to infinity) the random bipartite graph $G(n,n,p)$ has property $(\varepsilon, \varepsilon, s(n))$.
\end{cor}
In other words, almost all equi-bipartite graphs exhibit partial unimodality of the independent set sequence. (See \cite{Galvin4} for further results in this direction.)

\medskip

The bound in (\ref{eq-general_ub_on_part}) is valid for {\em all} $d$-regular graphs, and we expect that it, as well as both the upper and lower bounds on $i_t(G)$ ((\ref{eq-fromCGT}) and (\ref{eq-sillylb})), can be significantly improved if further structural conditions are put on $G$. In \cite{Galvin3} the programme of improving (\ref{eq-general_ub_on_part}) is carried out in the case when $G$ is the discrete hypercube $Q_d$. This is the
graph on vertex set $V=\{0,1\}^d$ with two strings adjacent if they
differ on exactly one coordinate. It is a $d$-regular bipartite
graph with bipartition classes $\cE$ and $\cO$, where $\cE$ is the
set of vertices with an even number of $1$'s. Note that
$|\cE|=|\cO|=\alpha(Q_d)=2^{d-1}$. The following bound is obtained in \cite{Galvin3}. (All asymptotic statements in what follows are as $d \rightarrow \infty$.)
\begin{thm} \label{thm-Zdest}
There is a $c>0$ and a function $f(d)\rightarrow 0$ such that for $\lambda > \frac{c\log d}{d^{1/3}}$,
$$
P(Q_d, \lambda) = 2(1+\lambda)^{2^{d-1}}\exp\left\{\frac{\lambda}{2}\left(\frac{2}{1+\lambda}\right)^d(1+f(d))\right\}.
$$
\end{thm}
This generalizes work of Korshunov and Sapozhenko \cite{KorshunovSapozhenko}, who had shown
$$
P(Q_d, 1) = (2\sqrt{e}+o(1))2^{2^{d-1}}.
$$

Using Theorem \ref{thm-Zdest} (or rather, using two of the intermediate inequalities that ultimately lead to the theorem), we can significantly improve the bounds on $i_t(Q_d)$ given by (\ref{eq-fromCGT}) and (\ref{eq-sillylb}), obtaining optimal (asymptotically matching) bounds for a large range of values of $t$.
\begin{thm} \label{thm-fixed_size}
There is a constant $c>0$ and a function $f(d)\rightarrow 0$ such that if $t=t(d)$ eventually (i.e., for all but finitely many $d$) satisfies
\begin{equation} \label{range-t4}
2^{d-1}\left(\frac{c\log d}{d^{1/3}}\right) \leq t \leq 2^{d-1}\left(1-\frac{1}{\sqrt{2}}+\frac{2\log d}{d}\right)
\end{equation}
then
$$
i_t(Q_d) = {2^{d-1} \choose t} \exp\left\{t \left(1-\frac{t}{2^{d-1}}\right)^{d-1}(1+f(d))\right\}.
$$
If $t$ eventually satisfies
\begin{equation} \label{range-t123}
2^{d-1}\left(1-\frac{1}{\sqrt{2}}+\frac{2\log d}{d}\right) \leq t \leq 2^{d-1}
\end{equation}
then
$$
i_t(Q_d) \sim 2{2^{d-1} \choose t} \exp\left\{t \left(1-\frac{t}{2^{d-1}}\right)^{d-1}\right\}.
$$
\end{thm}
\noindent After some groundwork in Section \ref{sec-NT}, we prove Theorem \ref{thm-fixed_size} in Section \ref{sec-proof_fixed}.

\medskip

A corollary of Theorem \ref{thm-fixed_size} is that the quantity $i_t(Q_d)$ undergoes a transition around $t=2^{d-2}$, in a window of width $O(1/d)$. This is analogous to \cite[Corollary 1.2]{Galvin3}.
\begin{cor} \label{cor-transition}
If $t=t(d) = 2^{d-1}\left(\frac{1}{2} + \frac{g(d)}{d}\right)$ then
$$
\frac{i_t(Q_d)}{2{2^{d-1} \choose t}} \rightarrow \left\{
\begin{array}{ll}
\infty & \mbox{if $g(d) \rightarrow -\infty$} \\
\exp\left\{e^{-2k}/2\right\} & \mbox{if $g(d) \rightarrow k$, a constant} \\
1 & \mbox{if $g(d) \rightarrow +\infty$}
\end{array}
\right.
$$
as $d \rightarrow \infty$.
\end{cor}

\medskip

We now turn to considering the unimodality of $P(Q_d,\lambda)$. By direct counting this polynomial may be shown to be unimodal for all $d \leq 5$ (although for $3 \leq d \leq 5$ it has some non-real roots).
For larger $d$, a corollary of Theorem \ref{thm-step_unimodality} is that for all $\varepsilon > 0$ there is a constant $C=C(\varepsilon)>0$ such that $Q_d$ has property $(0,\varepsilon,C2^d/d)$.
By examining the error terms in the asymptotic estimates of $i_t(Q_d)$ provided by Theorem \ref{thm-fixed_size}, we are able to obtain much stronger partial unimodality result for $P(Q_d, \lambda)$, and in particular extend Levit and Mandrescu's observation (\ref{seq-LM}) (in the particular case of $G=Q_d$) to a wider range of coefficients.
\begin{thm} \label{thm-unimodality}
For all sufficiently large $d$ it holds that
$$
i_p(Q_d) < i_{p+1}(Q_d) < \ldots < i_{2^{d-2}-15d^2}(Q_d)
$$
where $p=[(1-1/\sqrt{2}+2\log d/d)2^{d-1}]$, and
$$
i_{2^{d-2}+5d^4}(Q_d) > i_{2^{d-2}+5d^4+1}(Q_d) > \ldots > i_{2^{d-1}}(Q_d).
$$
\end{thm}
We give the proof in Section \ref{sec-proof_uni}. We have not made an attempt to optimize the coefficients of $d^2$ and $d^4$ here. The proof strategy is to show that for $p \leq t \leq t+1 \leq 2^{d-2}-15d^2$ our upper bound on $i_t(Q_d)$ from Theorem \ref{thm-fixed_size} is eventually smaller than our lower bound on $i_{t+1}(Q_d)$ (with a similar approach for $2^{d-2}+5d^4 \leq t \leq t+1 \leq 2^{d-1}$). Below $p$ the upper and lower bounds on $i_t(Q_d)$ provided by Theorem \ref{thm-fixed_size} are too far apart to be of any use, and in the neighbourhood of $t=2^{d-2}$ the sequence $\{i_t(Q_d)\}$ seems to be too flat for the present approach to be helpful.


\section{Proofs of Theorems \ref{thm-step_unimodality} and \ref{thm-step_unimodality-gen}} \label{sec-perc}

We begin with the proof of Theorem \ref{thm-step_unimodality}. Let $j$ and $\ell$ satisfy $0 \leq j \leq \ell \leq (1-\varepsilon)n/2$. By (\ref{eq-fromCGT}) and (\ref{eq-sillylb}), a sufficient condition for $i_\ell(G) > i_j(G)$ is
\begin{equation} \label{suff}
H\left(\frac{\ell}{n}\right)-H\left(\frac{j}{n}\right) > \frac{1}{d} + \frac{\log 2n}{2n}.
\end{equation}
By the mean value theorem,
$$
H\left(\frac{\ell}{n}\right)-H\left(\frac{j}{n}\right) = \left(\frac{\ell}{n} - \frac{j}{n}\right)H'(\xi)
$$
for some $\xi \in (\ell/n,j/n)$. Since the minimum of $H'(x)$ on $[0,(1-\varepsilon)/2]$ is achieved at $x=(1-\varepsilon)/2$ and is a positive  constant depending on $\varepsilon$, we get that a sufficient condition for $i_\ell(G) > i_j(G)$ is
$$
\ell-j \geq C(\varepsilon)\left(\frac{n}{d}+\log 2n\right).
$$
For $j$ and $\ell$ satisfying $(1+\varepsilon)n/2 \leq \ell \leq j \leq n$
the proof is almost identical,
once we observe that $H(x)$ is symmetric around $x=1/2$.

\medskip

In order to prove Theorem \ref{thm-step_unimodality-gen}, we need an analog of (\ref{eq-fromCGT}) for the family of $G$'s under consideration in that theorem.
The following result is a special case of a general result on graph homomorphisms from \cite[Section 3]{EngbersGalvin}. For $\lambda_1, \lambda_2 > 1$,
\begin{equation} \label{gen-funct}
\sum_{t=0}^{2n} i_t(G)\lambda_1^t \lambda_2^{2n-t} \leq (\lambda_1\lambda_2+\lambda_2^2)^n C(\lambda_1,\lambda_2)^{nh(G,d)}
\end{equation}
where $C(\lambda_1,\lambda_2)$ may be taken to be $\max \left\{256, \lambda_1+\lambda_2 \right\}$.
It follows that for all fixed $\lambda > 0$,
\begin{equation} \label{eq-fromEG}
P(G,\lambda) \leq (1+\lambda)^n C(\lambda)^{nh(G,d)}
\end{equation}
where $C(\lambda) > 0$ goes to infinity both as $\lambda$ goes to $0$ and as $\lambda$ goes to infinity. Indeed, if $\lambda \geq 1$ then we may take $\lambda_1=2\lambda$ and $\lambda_2=2$ in (\ref{gen-funct}) to get
$$
P(G,\lambda) = \lambda_2^{-2n} \sum_{t=0}^{2n} i_t(G)\left(\lambda_1/\lambda_2\right)^t \leq (1+\lambda)^n \max\left\{256, 2(1+\lambda)\right\}^{nh(G,d)}
$$
while if $\lambda < 1$ then we may take $\lambda_1=2$ and $\lambda_2=2/\lambda$ to get
$$
P(G,\lambda) \leq (1+\lambda)^n \max\left\{256, 2(1+1/\lambda)\right\}^{nh(G,d)},
$$
showing that we may take $C(\lambda)$ to be $2(1+\lambda)$ for $\lambda \geq 127$, to be $256$ for $1/127 \leq \lambda \leq 127$, and to be $2(1+1/\lambda)$ for $\lambda \leq 127$.

Using (\ref{eq-fromEG}) in place of (\ref{eq-general_ub_on_part}), we may reproduce the derivation of (\ref{eq-fromCGT}) to obtain the following upper bound on $i_t(G)$ for $G$ satisfying the conditions of Theorem \ref{thm-step_unimodality-gen}:
$$
i_t(G) \leq \exp_2\left\{H\left(\frac{t}{n}\right)n + C(t,n)nh(G,d)\right\}
$$
where $C(t,n)>0$ may be taken as follows:
$$
C(t,n) = \left\{
\begin{array}{ll}
\log \left(\frac{2n}{n-t}\right) & \mbox{if $\frac{127n}{128} \leq t$} \\
8 & \mbox{if $\frac{n}{128} \leq t \leq \frac{127n}{128}$} \\
\log \left(\frac{2n}{t}\right) & \mbox{if $t \leq \frac{n}{128}$}.
\end{array}
\right.
$$
It follows that for $\varepsilon > 0$ there is $C(\varepsilon) > 0$ such that in the range $t \in [\varepsilon n, (1-\varepsilon)n]$ we have
\begin{equation} \label{eq-gen-up-bound-perc}
i_t(G) \leq \exp_2\left\{H\left(\frac{t}{n}\right)n + C(\varepsilon)nh(G,d)\right\}.
\end{equation}
Noting that (\ref{eq-sillylb}) is still valid in the setting of Theorem \ref{thm-step_unimodality-gen}, it follows from (\ref{eq-gen-up-bound-perc}) (just as (\ref{suff}) followed from (\ref{eq-fromCGT})) that for all $\varepsilon > 0$, if $j$ and $\ell$ satisfy $\varepsilon n \leq j \leq \ell \leq (1-\varepsilon)n/2$, then a sufficient condition for $i_\ell(G) > i_j(G)$ is
$$
H\left(\frac{\ell}{n}\right)-H\left(\frac{j}{n}\right) > C(\varepsilon)h(G,d) + \frac{\log 2n}{2n}.
$$
The proof of Theorem \ref{thm-step_unimodality-gen} now goes through exactly as the proof of Theorem \ref{thm-step_unimodality}.

\section{Proofs of Theorems \ref{thm-fixed_size} and \ref{thm-unimodality}} \label{sec-cube}

\subsection{Preliminaries} \label{sec-NT}

We begin with some notation. For $A \subseteq V~(=\{0,1\}^d)$ write $N(A)$
for the set of vertices outside $A$ that are neighbours of a vertex
in $A$ and set
$$
[A]=\{v \in V:N(\{v\})\subseteq N(A)\}.
$$
Note that if $A$ is an independent set then $A \subseteq [A]$. Say that $A \subseteq \cE$ (or $\cO$) is {\em small} if $|[A]| \leq 2^{d-2}$ and {\em $2$-linked} if $A \cup N(A)$ induces a connected subgraph of $Q_d$.
Any $A$ can be decomposed into its maximal $2$-linked subsets; we refer to these as
the {\em $2$-components} of $A$, and write $k(A)$ for the number of $2$-components of $A$ and ${\rm cl}(A)$ for the size of the largest $2$-component of $A$. Finally, for all $\lambda > 0$ and $a, g \geq 0$ we define a function $F_\lambda(a,g)$ by
$$
F_\lambda(a,g) = \lambda^a(1+\lambda)^{-g}.
$$

There are two bounds from \cite{Galvin3} (intermediate steps in the derivation of Theorem \ref{thm-Zdest}) that we will make use of.
\begin{lemma} \label{lem-from_G_threshold}
There is a constant $c>0$ such that for all $\lambda > \frac{c\log d}{d^{1/3}}$
we have
\begin{equation} \label{eq-inner_sum_upper}
\sum_{A \subseteq \cE~{\rm small}} F_\lambda(|A|,|N(A)|) \leq \exp\left\{\frac{\lambda}{2}\left(\frac{2}{1+\lambda}\right)^d + \frac{d^2 \lambda^2 (1+\lambda)^2 2^d}{(1+\lambda)^{2d}}\right\}
\end{equation}
as well as, for fixed $k \geq 1$,
\begin{equation} \label{eq-inner_sum_upper_spec}
\sum_{A \subseteq \cE~{\rm small},~2{\rm-linked},~|A|\geq k} F_\lambda(a,g) \leq e^{k-1}d^{2k-2}2^d F_\lambda(k,kd-2k(k-1)).
\end{equation}
\end{lemma}

\medskip

\noindent {\em Proof}: The first inequality is \cite[Equation (23)]{Galvin3} and the second is \cite[Corollary 3.11]{Galvin3}. \qed


\subsection{Proof of Theorem \ref{thm-fixed_size}}
\label{sec-proof_fixed}

We assume throughout that
$$
t \geq \left(\frac{c\log d}{d^{1/3}}\right) 2^{d-1}
$$
with $c$ the same as the constant appearing in the range of validity of (\ref{eq-inner_sum_upper}) and (\ref{eq-inner_sum_upper_spec}). We also, where necessary, assume that $d$ is large enough to support our assertions. All asymptotic statements in what follows will be as $d \rightarrow \infty$.

\subsubsection{Lower bounds on $i_t(Q_d)$}

For the lower bounds in Theorem \ref{thm-fixed_size}, we may assume $t \leq \frac{3}{4}2^{d-1}$ (any constant greater than $1/2$ would do in place of $3/4$ here) since for $t \geq \frac{3}{4} 2^{d-1}$, $t\left(1-\frac{t}{2^{d-1}}\right)^{d-1} = o(1)$ and so the bound
$$
i_t(Q_d) \geq (2-o(1)){2^{d-1} \choose t}\exp\left\{t\left(1-\frac{t}{2^{d-1}}\right)^{d-1}\right\}
$$
is trivial. For the remaining range of $t$, let $f=f(t,d)$ be defined by
\begin{equation} \label{f-prop}
f = \max\left\{d,\, 5^7 e\gld\right\}.
\end{equation}
Both lower bounds in Theorem \ref{thm-fixed_size} will be based on the following inequality.
\begin{lemma} \label{lem-using-sym}
$$
i_t(Q_d) \geq 2\sum_{A \subseteq \cE, ~{\rm cl}(A) \leq 1, ~|A|\leq f} {2^{d-1}-d|A| \choose t-|A|}.
$$
\end{lemma}

\medskip

\noindent {\em Proof}: We get a lower bound on $i_t(Q_d)$ by first choosing a set of vertices of size no more than $f$, no two of which share a common neighbour, to be the intersection of the independent set with one of the two partition classes, and then extending it by choosing an arbitrary subset of the other class of appropriate size (note that for $A$ with ${\rm cl}(A)=1$, $|N(A)|=d|A|$). This gives
\begin{eqnarray*}
i_t(Q_d) & \geq & \sum_{A \subseteq \cE~\mbox{or}~A \subseteq \cO, ~{\rm cl}(A) \leq 1, ~|A|\leq f} {2^{d-1}-|N(A)| \choose t-|A|} \\
& = & 2\sum_{A \subseteq \cE, ~{\rm cl}(A) \leq 1, ~|A|\leq f} {2^{d-1}-d|A| \choose t-|A|}.
\end{eqnarray*}
The equality above is valid since $f < t/2$ and so there is no overlap between the independent sets of size $t$ which intersect $\cE$ in no more than $f$ vertices and those which intersect $\cO$ in no more than $f$ vertices. \qed

\medskip

For all $a$ and $g$ with $g \geq a$, and $0 \leq t \leq 2^{d-1}$, we have the identity
\begin{eqnarray}
{2^{d-1}-g \choose t-a} & = & {2^{d-1} \choose t} \frac{t^a(2^{d-1}-t)^{g-a}}{\left(2^{d-1}\right)^g} E(a,g) \nonumber \\
& = & F_{\lambda(t)}(a,g){2^{d-1} \choose t} E(a,g) \label{eq-introducing_gl}
\end{eqnarray}
where $\lambda(t) = t/(2^{d-1}-t)$ and
$$
E(a,g) = \frac{\prod_{i=0}^{a-1} \left(1-\frac{i}{t}\right)\prod_{i=0}^{g-a-1} \left(1-\frac{i}{2^{d-1}-t}\right)}
{\prod_{i=0}^{g-1} \left(1-\frac{i}{2^{d-1}}\right)}.
$$
Note that
\begin{equation} \label{int44}
F_{\lambda(t)}(k,dk) = \left(\frac{t}{2^{d-1}}\left(1-\frac{t}{2^{d-1}}\right)^{d-1}\right)^k,
\end{equation}
and that for fixed $d$ and $k$, the quantity on the right-hand side in (\ref{int44}) is decreasing in $t$ for the range of $t$ that we are considering (indeed, it is decreasing for $t \geq (2^{d-1}-1)/(d-1)$), a fact that we will use repeatedly in the calculations that follow, usually without comment.
For those $A$ contributing to the sum in Lemma \ref{lem-using-sym},
\begin{eqnarray}
E(|A|,d|A|) & \geq & \prod_{i=0}^{f-1} \left(1-\frac{i}{t}\right) \prod_{i=0}^{(d-1)f-1} \left(1-\frac{i}{2^{d-1}-t}\right) \nonumber \\
& \geq & \exp\left\{-\frac{f^2}{t}-\frac{d^2f^2}{2^{d-1}-t}\right\} \nonumber \\
& \geq & \exp\left\{-\frac{2f^2}{t}\right\}. \label{int1}
\end{eqnarray}
In the second inequality we use $1-x \geq e^{-2x}$ for $0 \leq x \leq 1/2$. The use is valid since both $f \leq t/2$ and $f \leq (2^{d-1}-t)/(2d)$ hold.

The number of ways of choosing $A \subseteq \cE$ with ${\rm cl}(A) \leq 1$ and $|A| = k \leq f$ is at least
\begin{eqnarray}
\frac{\prod_{i=0}^{k-1}\left(2^{d-1}-id^2\right)}{k!} & \geq & \frac{(2^{d-1})^k}{k!} \exp\left\{-\frac{f^2 d^2}{2^{d-1}}\right\} \nonumber \\
& \geq & \frac{(2^{d-1})^k}{k!} \exp\left\{-\frac{f^2}{t}\right\} \label{int2}
\end{eqnarray}
since each choice of vertex in $A$ eliminates from consideration at most $d^2$ other vertices. By (\ref{eq-introducing_gl}) each such $A$ contributes
\begin{equation} \label{int3}
F_{\lambda(t)}(k,dk){2^{d-1} \choose t} E(|A|, d|A|)
\end{equation}
to the sum in Lemma \ref{lem-using-sym}. Combining (\ref{int44}), (\ref{int1}), (\ref{int2}) and (\ref{int3}) with Lemma \ref{lem-using-sym} we get
\begin{equation} \label{int4}
i_t(Q_d) \geq 2 {2^{d-1} \choose t} \exp\left\{-\frac{3f^2}{t}\right\}  \sum_{k \leq f} \frac{1}{k!} \left(t\left(1-\frac{t}{2^{d-1}}\right)^{d-1}\right)^k.
\end{equation}
Using $k! \geq (k/e)^k$ and the lower bound on $f$ from (\ref{f-prop}) we have
$$
\sum_{k \leq f} \frac{1}{k!}\left(\gld\right)^k \geq E_0\exp\left\{\gld\right\}
$$
where
$$
E_0 = 1-2\left(\frac{e\gld}{f}\right)^f \exp\left\{-\gld\right\}
$$
and so (\ref{int4}) becomes
\begin{equation} \label{t_lower_bound}
i_t(Q_d) \geq 2{2^{d-1} \choose t} \exp\left\{t\left(1-\frac{t}{2^{d-1}}\right)^{d-1}\right\} E_1 
\end{equation}
where
$$
E_1  = \exp\left\{-\frac{3f^2}{t}\right\} E_0.
$$

The dominating term here is $\exp\{-3f^2/t\}$.
For $t$ satisfying (\ref{range-t4})
$$
E_1 \geq \exp\left\{-o\left(\gld\right)\right\}
$$
and for $t$ satisfying (\ref{range-t123}) $E_1 \geq 1-o(1)$, completing the lower bounds. For the purpose of proving Theorem \ref{thm-unimodality} we also note the following more precise bounds.
\begin{lemma} \label{lem-E1}
$$
E_1 \geq \left\{
\begin{array}{ll}
1-\frac{1}{d^5} & \mbox{if $t$ satisfies (\ref{range-t123})} \\
1-\frac{14d^2}{2^d} & \mbox{if $2^{d-1}\left(\frac{1}{2}-\frac{1}{d}\right) \leq t \leq \frac{3}{4}2^{d-1}$}.
\end{array}
\right.
$$
\end{lemma}

\subsubsection{Upper bounds on $i_t(Q_d)$}

We now turn to the upper bounds in Theorem \ref{thm-fixed_size}.
For any $I \in \cI(Q_d)$ we have $[I\cap \cE] \cap [I \cap \cO] = \emptyset$. Since $Q_d$ has a perfect matching, it follows that at least one of $[I \cap \cE], [I \cap \cO]$ is no larger than $2^{d-2}$, that is, that at least one of $I \cap \cE, I \cap \cO$ is small. This together with $\cE$-$\cO$ symmetry leads to the bound
\begin{equation} \label{fixed_ub}
i_t(Q_d) \leq 2 \sum_{A \subseteq \cE ~{\rm small}} {2^{d-1}-|N(A)| \choose t-|A|}.
\end{equation}
We will split the sum in (\ref{fixed_ub}) into three cases. Say that small $A \subseteq \cE$ is {\em of type I} if $|A|\leq f$, {\em of type II} if $|A| > f$ and ${\rm cl}(A) \leq 5$, and {\em of type III} if $|A| > f$ and ${\rm cl}(A) \geq 6$ (where $f$ is as defined in (\ref{f-prop})).

\medskip

\noindent {\bf Case 1 - Type I $A$'s}:
We first consider the contribution to the sum in (\ref{fixed_ub}) from $A$ of type I. For these $A$ we have
$$
E(|A|,|N(A)|) \leq \exp\left\{\frac{d^2f^2}{2^{d-1}}\right\}
$$
(this is similar to the derivation of (\ref{int1}), except that in this case we are lower bounding the numerator of $E(|A|,|N(A)|)$).
Taking $\lambda = \lambda(t) = t/(2^{d-1}-t)$ in (\ref{eq-inner_sum_upper}) and combining with (\ref{eq-introducing_gl})
we find that the contribution to (\ref{fixed_ub}) from $A$ of type I is at most
\begin{equation} \label{eq-fixed_ub_part_1}
2{2^{d-1} \choose t}  \exp\left\{t\left(1-\frac{t}{2^{d-1}}\right)^{d-1} + \frac{d^2 t^2 2^d}{(2^{d-1}-t)^2}\left(1-\frac{t}{2^{d-1}}\right)^{2d-2} + \frac{d^2f^2}{2^{d-1}}\right\}.
\end{equation}

\medskip

\noindent {\bf Case 2 - Type II $A$'s}:
Next we consider the sum in (\ref{fixed_ub}) over $A$ of type II. We have
\begin{eqnarray*}
{2^{d-1} - g \choose t - a} & = & \lambda(t)^{a-t} {2^{d-1} - g \choose t - a}\lambda(t)^{t-a} \\
& \leq & \lambda(t)^{a-t} (1+\lambda(t))^{2^{d-1} - g} \\
& = & (1+\lambda(t))^{2^{d-1}}\lambda(t)^{-t} F_{\lambda(t)}(a,g) \\
& = & 2^{H\left(\frac{t}{2^{d-1}}\right)2^{d-1}} F_{\lambda(t)}(a,g).
\end{eqnarray*}
By Stirling's formula, (more precisely, by the fact that for all $n\geq 1$,
$$
2n^ne^{-n}\sqrt{n} \leq n! \leq 3n^ne^{-n}\sqrt{n}),
$$
this is at most $3 \times 2^{d/2}{2^{d-1} \choose t}F_{\lambda(t)}(a,g)$. It follows that the contribution to (\ref{fixed_ub}) from $A$ of type II is at most
\begin{equation} \label{int33}
6\left(2^{d/2}\right){2^{d-1} \choose t} \sum_{A \subseteq \cE ~{\rm small},~{\rm cl}(A) \leq 5,~k(A)\geq f/5} F_{\lambda(t)}(|A|,|N(A)|).
\end{equation}
To bound this sum, we make a number of observations.
\begin{itemize}
\item For each $k \geq f/5$ there are at most $2^{k(d-1)}/k!$ ways of choosing a fixed vertex in each of the $k$ $2$-components of $A$, and at most $5^k$ ways of assigning a size to each $2$-component.
\item For each $\ell=1,\ldots, 5$, the number of $2$-linked subsets of $\cE$ of size $\ell$ that include a fixed vertex is at most $(\ell-1)!(d^2)^{\ell-1}$. (Once $j$ vertices have been chosen, there are at most $jd^2$ choices for the $(j+1)$st.)
\item Each $A \subseteq \cE$ with $|A| = \ell \leq 5$ satisfies $|N(A)| \geq d\ell-2\ell(\ell-1)$. (Each vertex of $A$ has $d$ neighbours,
of which at least $d-2(|A|-1)$ must be unique to it, since a pair of vertices in $Q_d$ can have at most two common neighbours.)
\item The quantity $F_\lambda(a,g)$ is decreasing in $g$, and for each $\ell=1, \ldots, 5$, (and sufficiently large $d$)
$$
(\ell-1)!(d^2)^{\ell-1} F_{\lambda(t)}(\ell,d\ell-2\ell(\ell-1)) \leq F_{\lambda(t)}(1,d).
$$
\end{itemize}
All this together serves to bound the expression in (\ref{int33}) by
\begin{equation} \label{eq-fixed_ub_part_2}
6\left(2^{d/2}\right){2^{d-1} \choose t} \sum_{k \geq f/5} \frac{1}{k!}\left(5t\left(1-\frac{t}{2^{d-1}}\right)^{d-1}\right)^{k} \leq \frac{1}{3^f}{2^{d-1} \choose t},
\end{equation}
the inequality following from the choice of $f$ and the bound $k!\geq (k/e)^k$.

\medskip

\noindent {\bf Case 3 - Type III $A$'s}:
Finally we consider the sum in (\ref{fixed_ub}) over $A$ of type III. Beginning with the same steps as in the case of $A$ of type II, this is at most
\begin{equation} \label{int27}
6\left(2^{d/2}\right){2^{d-1} \choose t} \sum_{A \subseteq \cE~{\rm small}, ~{\rm cl}(A)\geq 6} F_{\lambda(t)}(|A|,|N(A)|).
\end{equation}
We now use a multiplicative property of $F$: for any $\lambda$, if $A'$ is a $2$-component of $A$ then
$$
F_{\lambda}(|A|,|N(A)|) = F_{\lambda}(|A'|,|N(A')|)F_{\lambda}(|A\setminus A'|,|N(A \setminus A')|)
$$
and so the sum in (\ref{int27}) is at most
$$
\left(\sum_{A \subseteq \cE~{\rm small}, ~2{\rm-linked},~ |A|\geq 6} F_{\lambda(t)}(|A|,|N(A)|)\right) \left(\sum_{A \subseteq \cE~{\rm small}} F_{\lambda(t)}(|A|,|N(A)|)\right) = S_1 S_2.
$$
By (\ref{eq-inner_sum_upper}) we have
\begin{equation} \label{e1}
S_2 \leq \exp\left\{\gld + \frac{d^2 t^2 2^d}{(2^{d-1}-t)^2}\left(1-\frac{t}{2^{d-1}}\right)^{2d-2}\right\}
\end{equation}
and by (\ref{eq-inner_sum_upper_spec}) we have
\begin{equation} \label{e2}
S_1 \leq e^5d^{10}2^d F_{\lambda(t)}(6,6d-60).
\end{equation}

\medskip

Having examined the sum in (\ref{fixed_ub}) in each of the three possible cases, we now combine (\ref{eq-fixed_ub_part_1}), (\ref{eq-fixed_ub_part_2}), (\ref{e1}) and (\ref{e2}) to find that
\begin{equation} \label{t_upper_bound}
i_t(Q_d) \leq 2{2^{d-1} \choose t}  \exp\left\{t\left(1-\frac{t}{2^{d-1}}\right)^{d-1}\right\}E_2
\end{equation}
where
\begin{eqnarray*}
E_2 & = & \exp\left\{\frac{d^2 t^2 2^d}{(2^{d-1}-t)^2}\left(1-\frac{t}{2^{d-1}}\right)^{2d-2} + \frac{d^2f^2}{2^{d-1}}\right\} + \frac{1}{3^f} +\\
&  & ~~3e^5d^{10}2^{3d/2} F_{\lambda(t)}(6,6d-60) \exp\left\{\frac{d^2 t^2 2^d}{(2^{d-1}-t)^2}\left(1-\frac{t}{2^{d-1}}\right)^{2d-2}\right\}.
\end{eqnarray*}
For $t$ satisfying (\ref{range-t4})
$$
E_2 \leq \exp\left\{o\left(\gld\right)\right\}
$$
and for $t$ satisfying (\ref{range-t123}) $E_2 \leq 1+o(1)$, completing the upper bounds. For the purpose of proving Theorem \ref{thm-unimodality}, we also note the following more precise bounds.
\begin{lemma} \label{lem-E2}
$$
E_2 \leq \left\{
\begin{array}{ll}
1+\frac{1}{d^3} & \mbox{if $t$ satisfies (\ref{range-t123})} \\
1+\frac{4d^4}{2^d} & \mbox{if $t \geq 2^{d-1}\left(\frac{1}{2}-\frac{1}{d}\right)$}.
\end{array}
\right.
$$
\end{lemma}

\subsection{Proof of Theorem \ref{thm-unimodality}}
\label{sec-proof_uni}

We will split the proof into four cases, according to various possible ranges for $t$, specifically
\begin{equation} \label{t0}
(3/4)2^{d-1} \leq t < 2^{d-1},
\end{equation}
\begin{equation} \label{t1}
2^{d-2}+5d^4 \leq t < (3/4)2^{d-1},
\end{equation}
\begin{equation} \label{t2}
2^{d-1}\left(1-\frac{1}{\sqrt{2}}+\frac{2\log d}{d}\right) \leq t < 2^{d-1}\left(\frac{1}{2}-\frac{1}{d}\right),
\end{equation}
and
\begin{equation} \label{t3}
2^{d-1}\left(\frac{1}{2}-\frac{1}{d}\right) \leq t < 2^{d-2}-15d^2.
\end{equation}
In each of the first two cases we will show that for sufficiently large $d$ we have $i_t(Q_d)>i_{t+1}(Q_d)$, while in each of the last two cases we will show (again for sufficiently large $d$) that $i_t(Q_d) < i_{t+1}(Q_d)$.

\medskip

\noindent {\bf Case 1 - $t$ satisfying (\ref{t0})}: That $i_t(Q_d)>i_{t+1}(Q_d)$ for $t$ in this range follows from (\ref{seq-LM}).

\medskip

\noindent {\bf Case 2 - $t$ satisfying (\ref{t1})}: By (\ref{t_lower_bound}), (\ref{t_upper_bound}) and Lemmas \ref{lem-E1} and \ref{lem-E2}, for $t$ in this range it is enough to show that
$$
\frac{2\left(1-\frac{14d^2}{2^d}\right){2^{d-1} \choose t} \exp\left\{t\left(1-\frac{t}{2^{d-1}}\right)^{d-1}\right\}}{2\left(1+\frac{4d^4}{2^d}\right){2^{d-1} \choose t+1}  \exp\left\{(t+1)\left(1-\frac{t+1}{2^{d-1}}\right)^{d-1}\right\}} > 1.
$$
Since $t\left(1-\frac{t}{2^{d-1}}\right)^{d-1}$ is decreasing in $t$, this inequality is implied by
$$
\frac{\left(1-\frac{14d^2}{2^d}\right)(t+1)}{\left(1+\frac{4d^4}{2^d}\right)(2^{d-1}-t)} > 1.
$$
For $t$ satisfying (\ref{t1}) this is in turn implied by
$$
\frac{\left(1-\frac{14d^2}{2^d}\right)(2^{d-2}+ 5d^4+1)}{\left(1+\frac{4d^4}{2^d}\right)(2^{d-2}-5d^4)} > 1,
$$
which holds for all sufficiently large $d$.

\medskip

\noindent {\bf Case 3 - $t$ satisfying (\ref{t2})}: In this range we wish to show $i_t(Q_d)<i_{t+1}(Q_d)$.
Again by (\ref{t_lower_bound}), (\ref{t_upper_bound}) and Lemmas \ref{lem-E1} and \ref{lem-E2} it is enough to show that
\begin{equation} \label{ets2}
\frac{2\left(1+\frac{1}{d^3}\right){2^{d-1} \choose t} \exp\left\{t\left(1-\frac{t}{2^{d-1}}\right)^{d-1}\right\}}{2\left(1-\frac{1}{d^5}\right){2^{d-1} \choose t+1}  \exp\left\{(t+1)\left(1-\frac{t+1}{2^{d-1}}\right)^{d-1}\right\}} < 1.
\end{equation}
Writing $h(a,b)$ for $a(1-b/2^{d-1})^{d-1}$ we have
\begin{eqnarray}
h(t,t) - h(t+1,t+1) & \leq & h(t,t) - h(t,t+1) \nonumber \\
& = & h(t,t)\left(1-\left(1-\frac{1}{2^{d-1}-t}\right)^{d-1}\right) \nonumber \\
& \leq & h(t,t)\left(\frac{2(d-1)}{2^{d-1}-t}\right) \label{int55}
\end{eqnarray}
with (\ref{int55}) valid for sufficiently large $d$ as long as $2^{d-1}-t \geq 2$. For $t$ satisfying (\ref{t2}) we therefore have
\begin{eqnarray*}
\exp\left\{h(t,t)-h(t+1,t+1)\right\} & \leq & \exp\left\{h(t,t)\left(\frac{2(d-1)}{2^{d-1}-t}\right)\right\} \\
& \leq & \exp\left\{h\left(2^{d-3},2^{d-3}\right) \left(\frac{2(d-1)}{2^{d-1}-2^{d-2}}\right)\right\} \\
& \leq & 1+.76^d
\end{eqnarray*}
(for large enough $d$) and so (\ref{ets2}) in this range is implied by
$$
\frac{\left(1+\frac{2}{d^3}\right)(t+1)}{\left(1-\frac{1}{d^5}\right)(2^{d-1}-t)} < 1.
$$
This is in turn implied by
$$
\frac{\left(1+\frac{2}{d^3}\right)(2^{d-2}-\frac{2^{d-1}}{d}+1)}{\left(1-\frac{1}{d^5}\right)(2^{d-2}+\frac{2^{d-1}}{d})} < 1,
$$
which holds for all sufficiently large $d$.

\medskip

\noindent {\bf Case 4 - $t$ satisfying (\ref{t3})}:
Again by (\ref{t_lower_bound}), (\ref{t_upper_bound}) and Lemmas \ref{lem-E1} and \ref{lem-E2} it is enough to show that
\begin{equation} \label{ets3}
\frac{2\left(1+\frac{4d^4}{2^d}\right){2^{d-1} \choose t} \exp\left\{t\left(1-\frac{t}{2^{d-1}}\right)^{d-1}\right\}}{2\left(1-\frac{14d^2}{2^d}\right){2^{d-1} \choose t+1}  \exp\left\{(t+1)\left(1-\frac{t+1}{2^{d-1}}\right)^{d-1}\right\}} < 1.
\end{equation}
For $t$ satisfying (\ref{t3}) we use (\ref{int55}) to obtain
$$
\exp\left\{h(t,t)-h(t+1,t+1)\right\} \leq 1+\frac{d^2}{2^d}
$$
and so (\ref{ets3}) is implied by
$$
\frac{\left(1+\frac{5d^4}{2^d}\right)(t+1)}{\left(1-\frac{14d^2}{2^d}\right)(2^{d-1}-t)} < 1.
$$
This in turn is implied by
$$
\frac{\left(1+\frac{5d^4}{2^d}\right)(2^{d-2}-15d^2+1)}{\left(1-\frac{14d^2}{2^d}\right)(2^{d-2}+15d^2)} < 1
$$
which holds for sufficiently large $d$.


\bigskip

{\em Acknowledgement}: Part of this work was carried out while the author was a participant in the programme on Combinatorics and Statistical Mechanics at the Isaac Newton Institute for Mathematical Sciences, University of Cambridge, in spring 2008. The author thanks the institute and the programme organizers for the support provided.

\end{document}